# A Minimal 7-Fold Rhombic Tiling

Theo P. Schaad

## Abstract:

A set of tiles for covering a surface is composed of three types of tiles. The base shape of each one of them is a diamond or rhombus, with acute angles of π/7, 2π/7, and 3π/7 radians (A, B, and C tiles). Such a set is called 7-fold, as the B tiles can be arranged in 7-fold rotational symmetry forming a 7-pointed star. A substitution scheme is discovered wherein a number of proto-tiles are arranged to fill larger similar tiles that have the same matching characteristics as the base tiles. This process can then be repeated to tile an arbitrarily large surface. 7-fold sets have been discovered before, but they require a large number of base tiles, at least 41 tiles or more. The set described here is thought to be the most minimal or smallest set wherein only 11 base tiles are needed to start the process.

Background:

Tile sets with n-fold symmetry have been studied extensively since the discovery of the Penrose 5-fold aperiodic pattern[1]. The 5-fold pattern is aesthetically pleasing because it incorporates the golden ratio, a number that appears in nature, architecture, and in the Fibonacci series. In 5-fold stars made with 5 rhombs (with unit side lengths), the long diagonal of the rhomb is the golden ratio. In similar fashion, the number 7, the 7-pointed star, and other heptagrams have significant religious and occult symbolism and have been widely used in heraldry and logos. The equivalent 7-fold magic number Φ (Phi) is the long diagonal of the B tile that can be tiled into a 7-pointed star. Finding 7-fold patterns could therefore be a worthwhile endeavor, especially if it incorporates 7-pointed stars. However, finding such patterns is far from trivial.

There are two known schemes to find large (infinite) tilings. The first one was pioneered by Penrose and involves substitution and matching rules to generate larger tiles that have the same characteristics as the base tiles. This process is called inflation and can be repeated as a recurring series. The other method is a multi-grid method developed by de Bruijn[2] that leads to a potentially infinite variety of solutions. P. Steinhardt[3] used a hepta-grid to find a particularly attractive tiling with 7-pointed stars. A reconstructed patch is shown in Fig. 1 to point out some general features about 7-fold tilings. 1. It is obvious that the three base tiles are not only sufficient to cover a large surface, all three are necessary. 2. In the Steinhardt patch are some straight lines that follow edges and diagonals of base tiles. This is promising for the substitution method, as one tries to place entire base tiles, or at least half-tiles along the straight edges of the inflated tile. 3. The Steinhardt tiling does not contain straight sequences that lend themselves to a substitution scheme.

An impressive 7-fold tiling with substitution and matching rules was discovered by J. Socolar[4,5]. The substitution requires 129 base tiles. Fig. 2 shows another substitution tiling that requires 353 base tiles discovered by the author. Socolar had noted that some combinations of tiles can easily be interchanged (are weakly matched). The advantage of large substitution schemes is that certain features, like 7-pointed stars, can be enhanced.

On the other hand, 7-fold tilings with a low inflation factor are much more constrained. An example was discovered by Chaim Goodman-Strauss[6]. A minor variation of his tiling is shown in Fig. 3 with added notches to specify the matching directions. It takes 41 base tiles to cover the three inflated tiles. The sides of the enlarged tiles follow the inflation sequence δ = 1 + Φ + 1, where Φ is the 7-fold magic number mentioned above. The inflation factor is δ = 3.802… It is noted that there are no 7-pointed stars in the Goodman-Strauss tiling. Fig. 4 is a patch of the tiling derived from Fig. 3 arranged as a rosette or Mandala in circular or rotational symmetry (a Mandala is a geometric configuration of symbols arranged by a division of a circle).

## A Minimal 7-fold Rhombic Tiling:

One goal of this study of 7-fold tilings was to find the most minimal or smallest set of base tiles that can be used to generate an infinite tessellation by substitution and matching rules. Such an example is shown in Fig. 5. Substitution rules set the number of tiles needed to tile the three larger tiles. In addition, matching rules are used to specify the orientation of the base tiles and how they can be adjoined with other tiles.

The substitution rule of this minimal tiling is as follows.

$A_1 = 2 A_0 + B_0$

$B_1 = A_0 + 2 B_0 + C_0$

$C_1 = B_0 + 3 C_0$

This can be conveniently written as

$$\begin{matrix} A_1 \\ B_1 \\ C_1 \end{matrix} = M \begin{matrix} A_0 \\ B_0 \\ C_0 \end{matrix}$$

Where M is the substitution matrix. In this case,

$$M = \begin{matrix} 2 & 1 & 0 \\ 1 & 2 & 1 \\ 0 & 1 & 3 \end{matrix}$$

If the same substitution can be repeated to the $n^{th}$ generation then

$$\begin{pmatrix} A_n \\ B_n \\ C_n \end{pmatrix} = M^n \begin{pmatrix} A_0 \\ B_0 \\ C_0 \end{pmatrix}$$

For instance, the substitution matrix of the second generation is

$$M^2 = \begin{pmatrix} 5 & 4 & 1 \\ 4 & 6 & 5 \\ 1 & 5 & 10 \end{pmatrix}$$

In this minimal tiling set, the sides of the inflated rhombs are enlarged by the inflation factor δ, which is the long diagonal of the A tile (denoted by Γ). The base tiles have sides of unit length.

Γ = 2 cos(π/14) (= 1.950…)

Hence the inflation factor of this minimal tiling is δ = 1.950…

Surprisingly, $M^2$ is the same substitution matrix as the Goodman-Strauss tiling of Fig. 4 which has an inflation factor of δ = Φ + 2 (= 3.802…). From this, one can derive one of many relationships between diagonals, areas, and area ratios:

$Γ^2$ = Φ + 2      (the long diagonal of the A tile is related to the 7-fold magic number Φ)

Figs. 5-8 show the first four generations of the minimal set. A variation is shown in Figs. 9-12, wherein the base tiles are taken from the first generation of the previous set. The $A_0$ tile ("shield") is composed of three and the $B_0$ ("skate") and $C_0$ tiles ("kite") are composed of four proto-rhombs. The name "kite" is in analogy to the kites and darts in the Penrose tiling, a moniker originally coined by John Horton Conway[7]. Another variation is shown in Figs. 14-16 nicknamed "minimal polka dot 7-fold tiling". It consists as before of only three proto-rhombs, but the coloring brings out a new set of tiles, two trapezoids, and one triangle.

## Magic and Impossible Substitutions:

Fig. 18 demonstrates that it is physically possible (with paper and scissors) to tile inflated tiles that are enlarged by the 7-fold magic number Φ. However, some additional cuts are necessary to create inflated rhombs. There are no known matching rules to make an infinite pattern. But the substitution matrix is perfectly valid:

$$M_Φ = \begin{pmatrix} 1 & 0 & 1 \\ 0 & 2 & 1 \\ 1 & 1 & 2 \end{pmatrix}$$

Φ is the long diagonal of the B tile and is the inflation factor of this "magic" substitution:

$\Phi = 2 \cos(\pi/7)$ (= 1.802...)

Surprisingly, the "magic" substitution matrix has a square root with integer value:

$$M_{\sqrt{\Phi}} = \begin{matrix} 0 & 1 & 0 \\ 1 & 0 & 1 \\ 0 & 1 & 1 \end{matrix}$$

Potentially, it should be possible to make inflated tiles with only 5 proto-tiles with an inflation factor of $\delta = \Phi^{1/2}$ (= 1.342...). But how can this be achieved? This is akin to asking, how can a circle be tiled with square tiles? The first substitution (top row of the matrix) is inflating the $A_0$ tile to the size of $B_0$. The ratio of the two areas is $B_0/A_0 = \Phi$, an irrational number. There is no way to tile it with only one proto-tile as it is impossible to tile a circle with a finite number of square tiles.

The "square root of $\Phi$" matrix leads to an interesting cubic equation. From the second row we get $B_1 = A_0 + C_0$. The ratio of the areas $B/A = \Phi$ (a summary of various relationships is given at the end of the paper). The area inflation is $\Phi$, i.e., $B_1 = \Phi B_0$, or $B_1 = \Phi^2 A_0$. From the second row, divided by $A_0$, $\Phi^2 = 1 + C_0/A_0$. The ratio of the areas $C/A = (C/B)(B/A) = \Phi(\Phi^2-2)$.

Rearranging all terms leads to the cubic equation:

$\Phi^3 - \Phi^2 - 2\Phi + 1 = 0$

This really defines the magic number $\Phi$ in the context of 7-fold tessellations. The equation has three real roots and $\Phi$ = 1.802... is one of them. This can be easily checked numerically. The equation is hard to solve, and I leave it as a challenge to the reader.

## The General Form of the 7-fold Substitution Matrix:

The $A_0$ tile could also be inflated to the size of the $C_0$ tile. Again, this is physically impossible (with paper and scissors) but makes perfect sense mathematically. The substitution matrix is:

$$M_C = \begin{matrix} 0 & 0 & 1 \\ 0 & 1 & 1 \\ 1 & 1 & 1 \end{matrix}$$

It has an inflation factor of $\delta = (\Phi^2-1)^{1/2}$ (= 1.499...). This can be extended to any integers a,b,c.

$$M_A = \begin{matrix} a & 0 & 0 \\ 0 & a & 0 \\ 0 & 0 & a \end{matrix} \quad M_B = \begin{matrix} 0 & b & 0 \\ b & 0 & b \\ 0 & b & b \end{matrix} \quad M_C = \begin{matrix} 0 & 0 & c \\ 0 & c & c \\ c & c & c \end{matrix}$$

And the combined substitution is M = $M_A$ + $M_B$ + $M_C$. Any 7-fold substitution matrix is therefore of the form

$$M = \begin{pmatrix} a & b & c \\ b & (a+c) & (b+c) \\ c & (b+c) & (a+b+c) \end{pmatrix}$$

It suffices to know the seed values (a,b,c) to tile $A_1$ to predict how many tiles are needed to tile $B_1$ and $C_1$. (a,b,c) can be any integer combination, but only a few will lead to a substitution and matching scheme that tiles an infinite area.

## The Magic Number Φ in an Infinite Series:

In a large 7-fold tiling, the number of tiles $N_A$, $N_B$, $N_C$ are in proportion of their areas, i.e.,

$N_B/N_A \sim B/A = Φ$ (large n)

Looking at just one of the inflated tiles in the nth generation, $A_n$, it is composed of

$A_n = n_A A_0 + n_B B_0 + n_C C_0$

The ratio of $n_B$ tiles to $n_A$ tiles approaches

$\lim(n \to \infty) n_B/n_A = Φ$

Any substitution matrix with an arbitrary ($a_0, b_0, c_0$) seed could lead to an infinite series that generates $n_A$ and $n_B$. Starting with the general substitution matrix M (above)

$A_1 = a_0 A_0 + b_0 B_0 + c_0 C_0$

Performing another inflation with M results in

$A_2 = b_0 A_0 + (a_0 + c_0) B_0 + (b_0 + c_0) C_0$

Which can be the start of another 7-fold tiling with

$A_2 = a_1 A_0 + b_1 B_0 + c_1 C_0$

And so on. This can be given as the following infinite series[8]:

Let ($a_0, b_0, c_0$) be any integer seed values

For n>0

$a_n = b_{n-1}$

$b_n = a_{n-1} + c_{n-1}$

$c_n = b_{n-1} + c_{n-1}$

then

$\lim(n \to \infty) b_n/a_n = \Phi$

A good starting seed is $a_0 = 5$, $b_0 = 9$, $c_0 = 11$ ($\Phi \sim a_0/b_0$), but any other seed will also converge.

## The 7-fold Magic Number Φ:

During this study of 7-folds, it became obvious that areas, area ratios, and diagonals of the three proto-rhombs are all related to Φ. Rhomb edges are normalized to unit length.

Areas:

| Tile | Area | Value | Related to Φ |
|---|---|---|---|
| A | $A = \sin(\pi/7)$ | 0.434 | $(\Phi+2)^{1/2}/(2\Phi^2-2)$ |
| B | $B = \sin(2\pi/7)$ | 0.782 | $\Phi(\Phi+2)^{1/2}/(2\Phi^2-2)$ |
| C | $C = \sin(3\pi/7)$ | 0.975 | $(\Phi+2)^{1/2}/2$ |

Area Ratios:

| Ratio | Value | Related to Φ |
|---|---|---|
| B/A | 1.802 | $\Phi$ |
| C/B | 1.247 | $\Phi^2-2$ |
| C/A | 2.247 | $\Phi^2-1$ |

Diagonals:

| Tile | Diag. | Trig. | Value | Related to Φ |
|---|---|---|---|---|
| A | Short | $\gamma = 2\sin(\pi/14)$ | 0.445 | $\Phi(\Phi^2-3)$ |
|   | Long | $\Gamma = 2\cos(\pi/14)$ | 1.950 | $(\Phi+2)^{1/2}$ |
| B | Short | $\phi = 2\sin(\pi/7)$ | 0.868 | $(\Phi+2)^{1/2}/(\Phi^2-1)$ |
|   | Long | $\Phi = 2\cos(\pi/7)$ | 1.802 | $\Phi$ |
| C | Short | $\psi = 2\sin(3\pi/14)$ | 1.247 | $\Phi^2-2$ |
|   | Long | $\Psi = 2\cos(3\pi/14)$ | 1.564 | $\Phi(\Phi+2)^{1/2}/(\Phi^2-1)$ |

## Summary:

In the following table, some of the 7-fold tilings mentioned in this paper are summarized and grouped by the inflation factor. The substitution matrix can be calculated from the seeds a,b,c where $A_1 = a A_0 + b B_0 + c C_0$ is the substitution of the smallest tile which determines the number of tiles needed for $B_1$ and $C_1$. The inflation factor δ is expressed in relation to the magic number Φ. The # Tiles is the total number of base tiles needed to tile the next generation. The label Pattern # is arbitrary and was roughly in the sequence of studying the particular pattern.

| Group | a,b,c | Inflation factor δ | = | # Tiles | Notes |
|---|---|---|---|---|---|
| A | 1,0,0 | 1 | 1 | 3 | Identity matrix (no inflation) |
| B | 0,1,0 | $\Phi^{½}$ | 1.342… | 5 | Inflate $A_0$ area to $B_0$ area |
| C | 0,0,1 | $(\Phi^2-1)^{½}$ | 1.499… | 6 | Inflate $A_0$ area to $C_0$ area |
| D | 1,0,1 | $\Phi$ | 1.802… | 9 | "Magic" inflation, shown in Fig. 18 |
| E | 2,1,0 | $(\Phi+2)^{½}$ | 1.950… | 11 | "Minimal" inflation, shown in Figs. 5-13<br>Binary Tilings[10] |
| F | 4,0,0 | 2 | 2 | 12 | Replace each tile with 4 of the same shape (a periodic pattern) |
| G | 2,2,1 | $\Phi+1$ | 2.802… | 22 | Pattern 1 shown in Fig. 17 |
| H | 5,4,1 | $\Phi+2$ | 3.802… | 41 | Goodman-Strauss 7-fold<br>Pattern 2 shown in Figs. 3,4 (similar) |
| I | 3,5,6 | $\Phi(\Phi+1)$ | 5.049… | 70 | Pattern 5 (not shown, could not find matching rules)<br>Madison[9] 7-fold |
| J | 6,9,11 | $\Phi(\Phi+2)$ | 6.851… | 129 | Socolar's 7-fold<br>Pattern 3 shown in Fig. 18 (different) |
| K | 9,16,20 | $\Phi^2(\Phi+1)$ | 9.098… | 227 | Pattern 4 (not shown, could not find matching rules) |
| L | 14,25,31 | $\Phi^3+2\Phi^2-1$ | 11.345… | 353 | Pattern 6 shown in Fig. 2 (a beauty with lots of 7-pointed stars) |
| M | a,b,c | $(a-c+b\Phi+c\Phi^2)^{½}$ | | 3a+5b+6c | General formula for 7-fold substitutions |

Epilogue:

This study of 7-fold tilings started on Valentine's Day 2020. Within weeks, the author was in quarantine for 3 months with his family due to the Coronavirus pandemic. It was not quite planned that way but cutting out hundreds of tiles and gluing them into patterns was a perfect pastime. The author expects that many other tiling enthusiasts have also tried their hands with 7-folds. In self-isolation, this work was independent, but the hope is to give credit where credit is due. I had started working on Pattern 1 and 2 in 1983 but did not get anywhere. I was surprised how difficult it was to find new patterns by substitution. But then again, every time I found a new 7-fold tessellation, I was amazed that it existed.

Sadly, the pandemic took its toll on John Horton Conway[7] who passed away in April. He made contributions to many branches of recreational mathematics and played an important role in the early days of the Penrose 5-fold tessellation and in solutions for Rubik's Cube.

Postscript:

A year after the first version of this paper was published, the author would like to acknowledge the feedback he received from Jim Millar[11] and Dieter Steemann[12]. Jim Millar pointed out that the minimal 7-fold tiling was known as a Binary Tiling published[10] by F. Lancon and L. Billard in 1993. The term "binary tiling" was originally used for binary Penrose tiles as a model for binary metal alloys. In the context of 7-fold rhombic tilings, the term ternary might be more appropriate, as there are three base tiles. Jim Millar posted 7-fold binary tilings on a website.

Dieter Steemann found a typo in the inflation factor of Group K of the summary table. It should be $\Phi^2(\Phi+1)$. He also shared his extensive notebook on 7-fold tilings that includes references and graphic computer code for Mathematica 11. Several tilings are listed in the Bielefeld Tiling Encyclopedia[14] including the aforementioned Harris 7-fold and binary tilings. In addition, there is another highly complex 7-fold rhombic tiling by A.E. Madison[9] and a recent one by Gregory Maloney. Other related tilings have base tiles that differ from rhombs, i.e., stars (S. Pautze), triangles (L. Danzer, S. Pautze, G. Maloney), and dihedral base tiles (J. Millar).

A somewhat different approach on a parquet tiling with 7-fold rotational symmetry was recently published by A. Wünsche[14] who shared his delight and frustration with the number "seven".

Finally, it can be shown that ternary rhombic 7-fold tilings (with three base tiles) can be written as an eigenmatrix equation with a characteristic eigenvalue that is the square of the inflation factor. The eigenvalue can be solved with a cubic differential equation, similar to the one shown in the paper. Each group in the summary table thus leads to a different cubic differential equation.

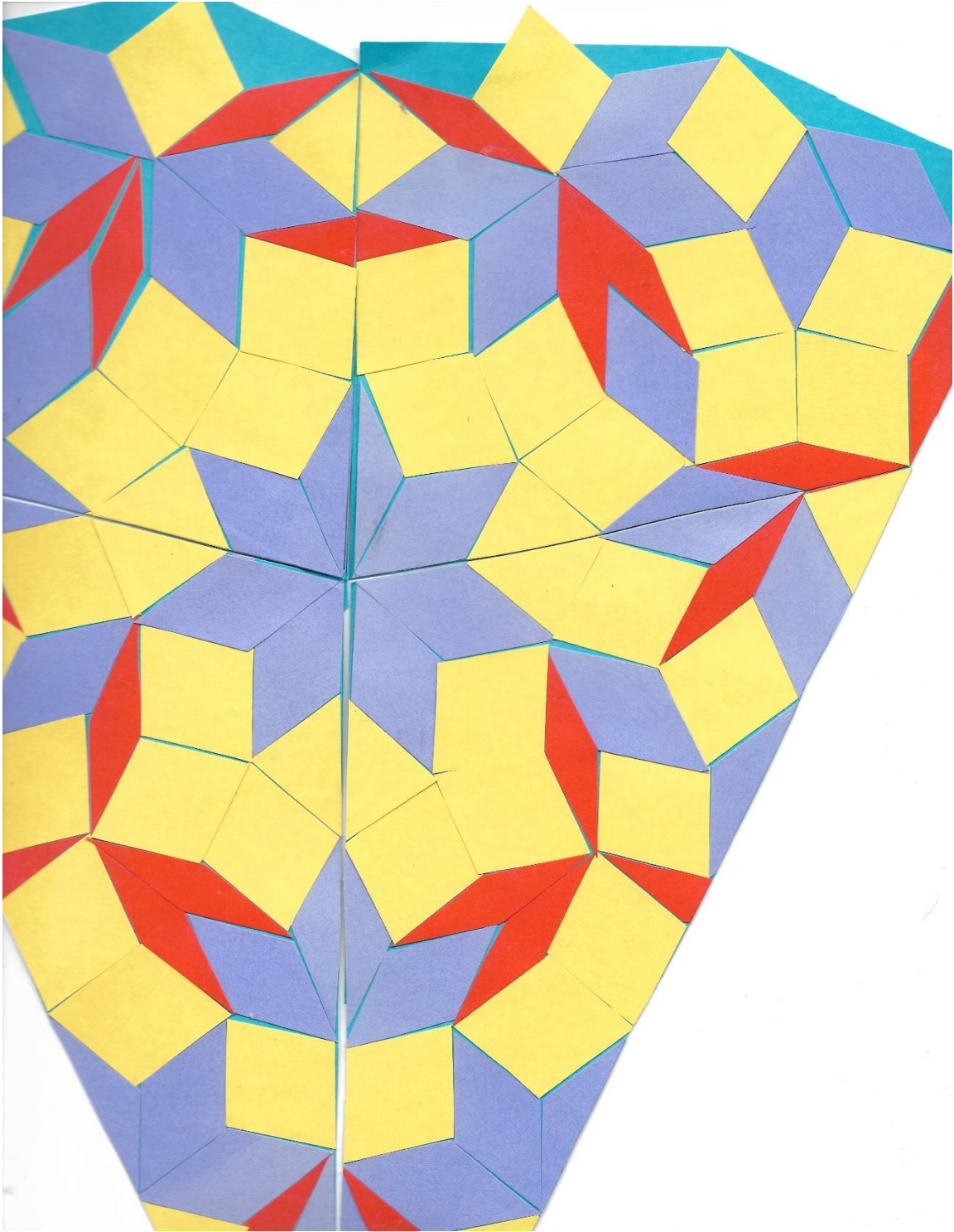

Fig 1: A patch of the Steinhardt 7-fold tiling constructed with the de Bruijn hepta-grid method.

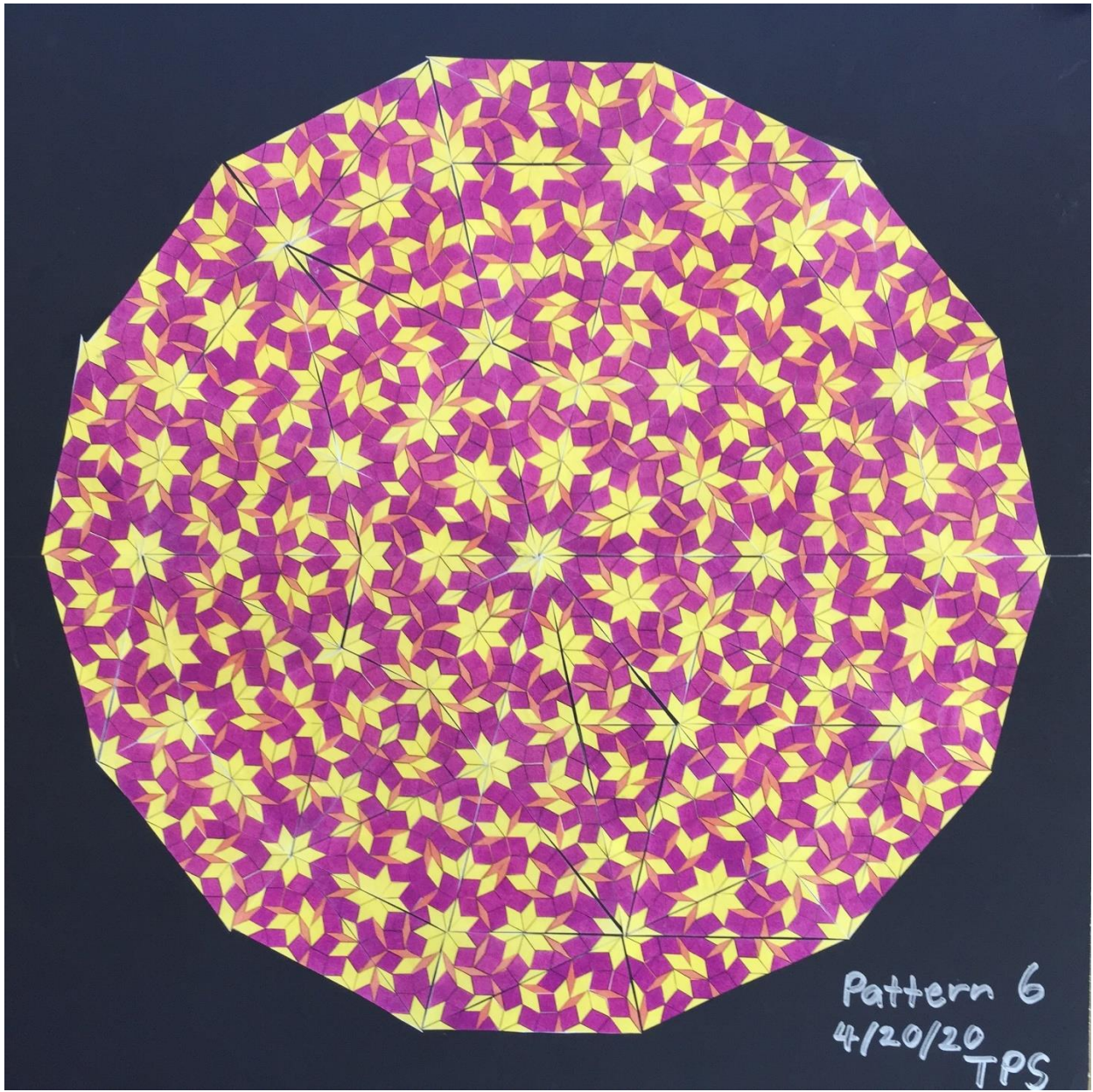

Fig 2: A central rosette of a 7-fold tiling constructed with a substitution and matching scheme. 353 base tiles are used to tile three similar larger tiles. The inflation sequence (sides of the inflated tiles) is δ = 1 + ψ + Φ + γ + 1 + Φ + ψ + 1 + Φ, where Φ is the long diagonal of the B tile, ψ is the short diagonal of the C tile, and γ is the short diagonal of the A tile. The inflation factor is δ = 11.345…

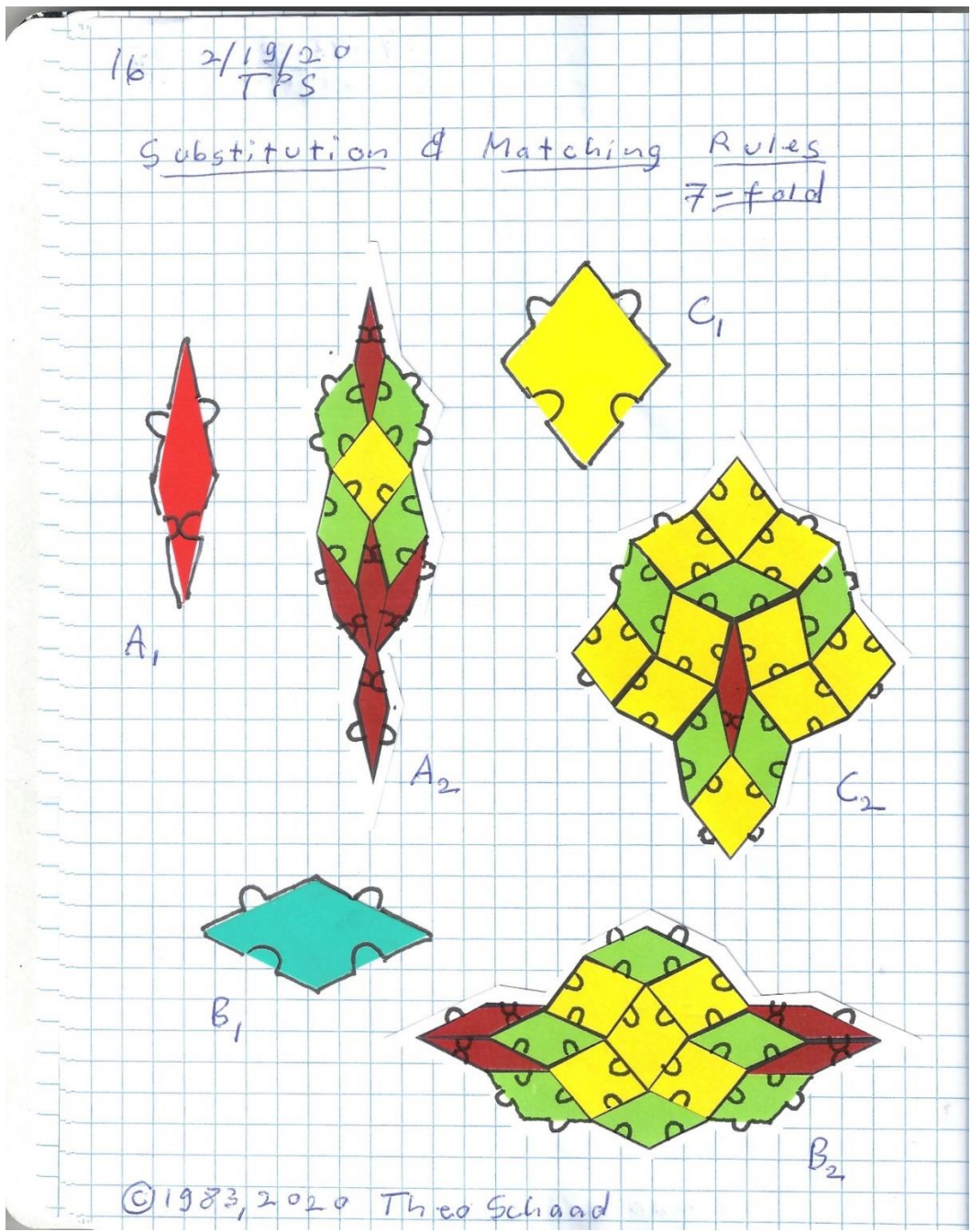

Fig 3: Substitution and matching rules for a 7-fold tiling. This pattern was found independently but is essentially the same as the one discovered by Chaim Goodman-Strauss. It takes 41 base tiles to make the three inflated tiles. In addition, notches were added to specify the matching directions. The inflation sequence (sides of the larger tiles) is δ = 1 + Φ + 1, and the inflation factor δ = 3.802...

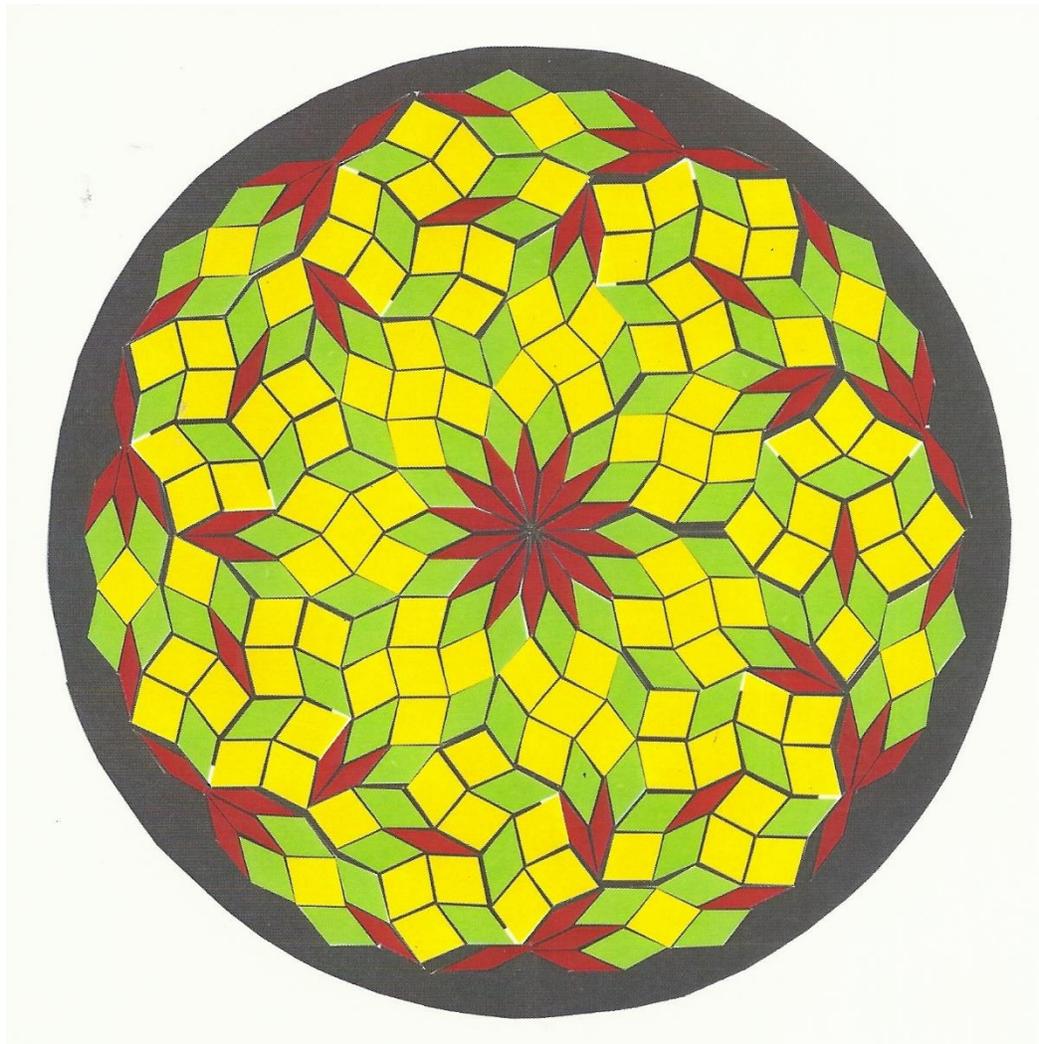

Fig 4: A central rosette of the 7-fold tiling of Fig. 3. There are no 7-pointed stars in this tiling.

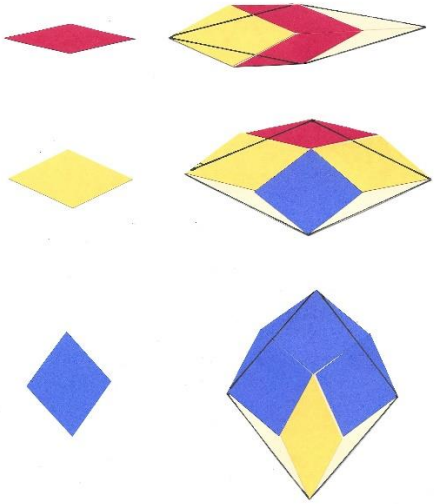

Fig. 5: Substitution rules of the minimal 7-fold rhombic tiling. The edges of the inflated tiles are now the length of the long diagonal of the A tile. The inflation factor is δ = 1.950...

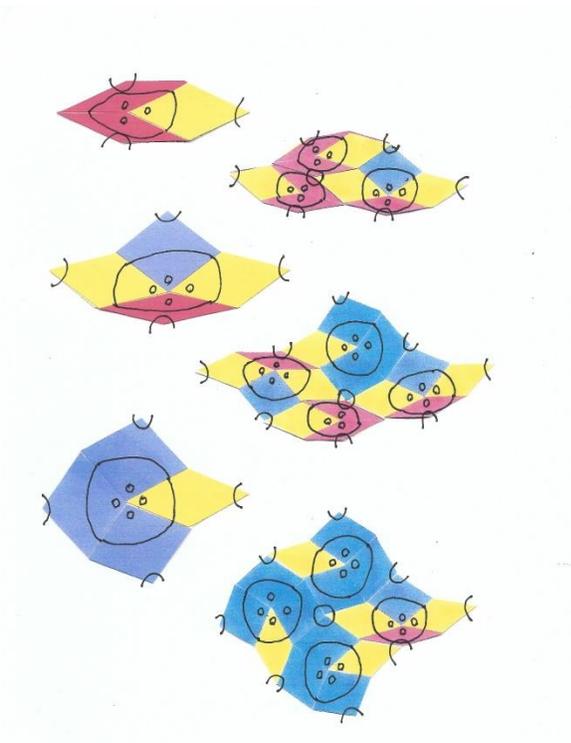

Fig. 6: Matching rules of the first and second generation of the minimal 7-fold tiling.

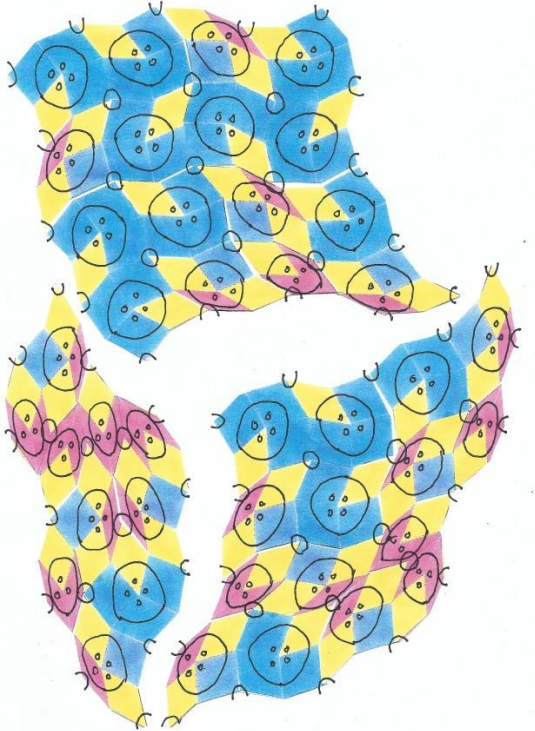

Fig. 7: The third generation of the minimal 7-fold tiling with substitution and matching rules.

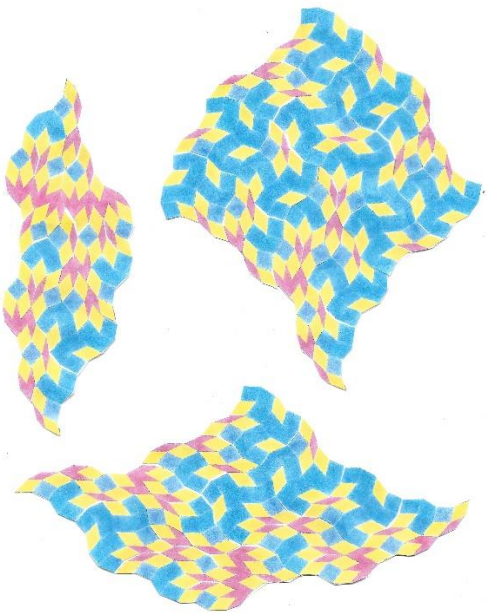

Fig. 8: The 4th generation of the minimal 7-fold rhombic tiling.

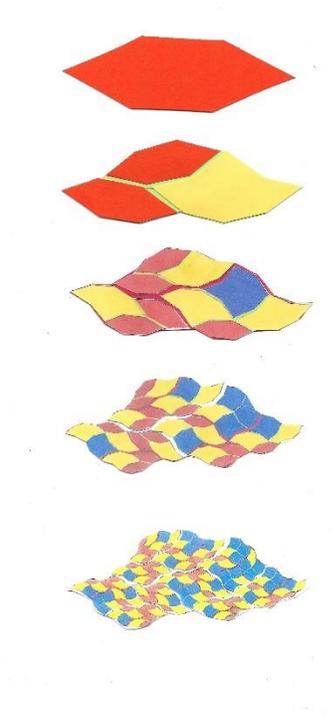

Fig. 9: A variation of the minimal 7-fold tiling. Here the base A tile ("shield") is taken from the first generation of Fig. 5 and is a composite of three smaller rhombs.

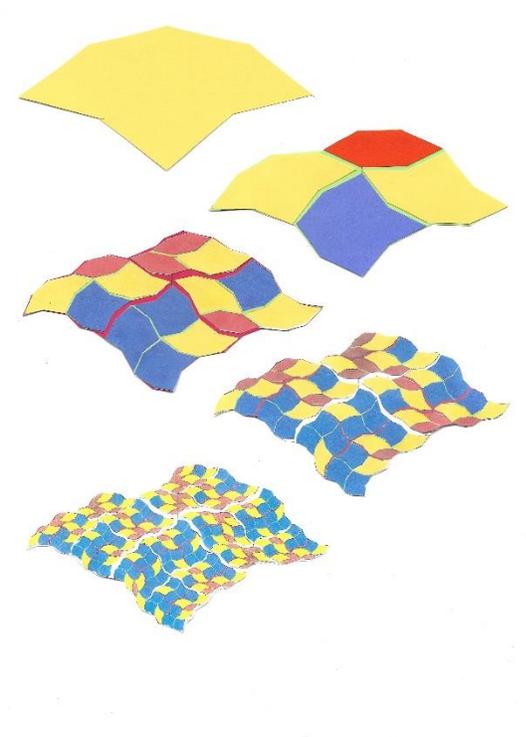

Fig. 10: Here the base B tile ("skate") is from the first generation of Fig. 5 and is itself a combination of four smaller rhombs.

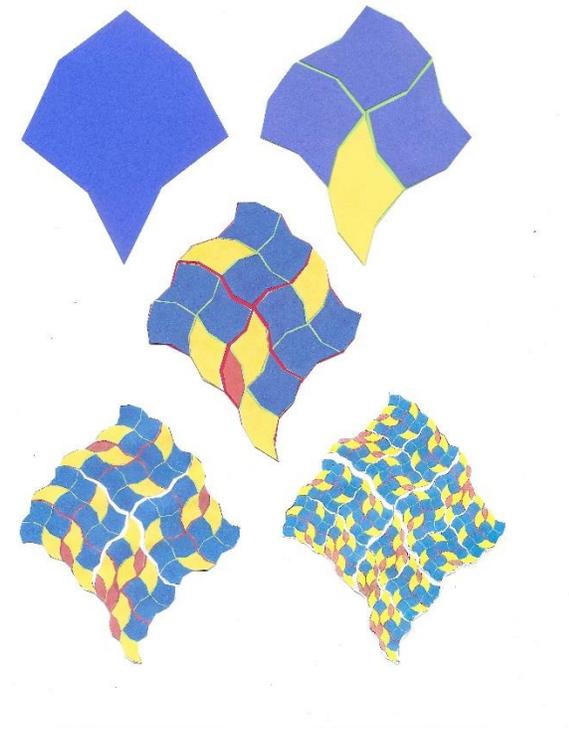

Fig. 11: The C tile ("kite") is taken from the first generation of Fig. 5 and is a combination of four smaller rhombs. Hence the pattern is a variation of the minimal 7-fold rhombic tiling.

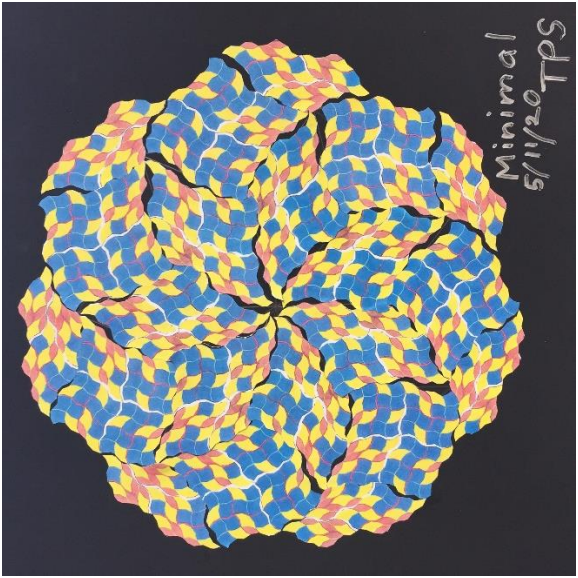

Fig. 12: A mandala of the minimal 7-fold tiling using shields, skates, and kites.

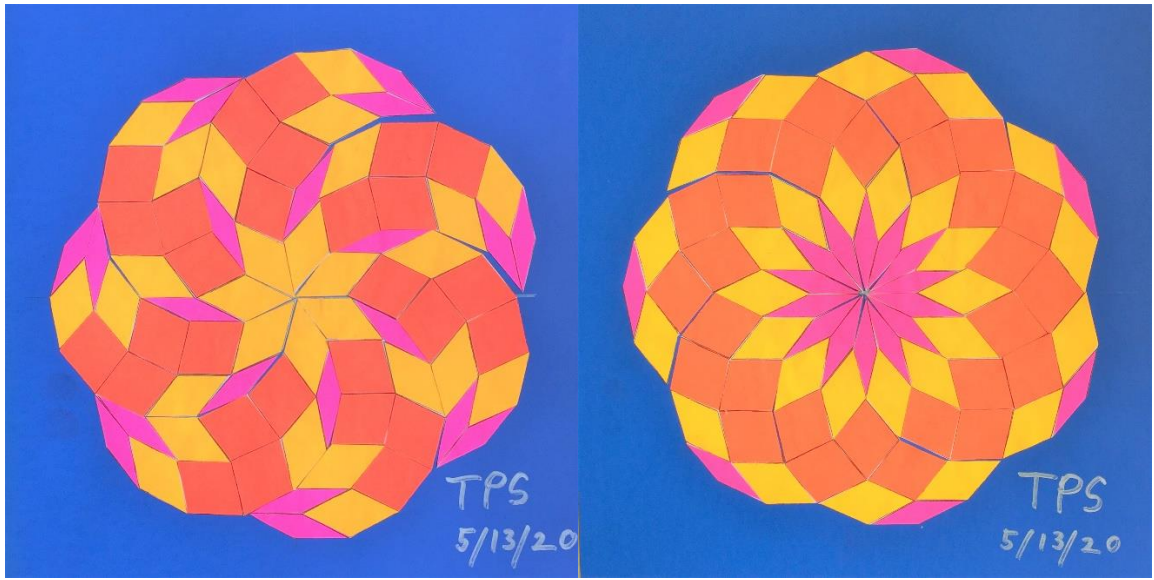

Fig. 13: Two mandalas of the minimal 7-fold tiling using the inflated tiles of the first generation in Fig. 5.

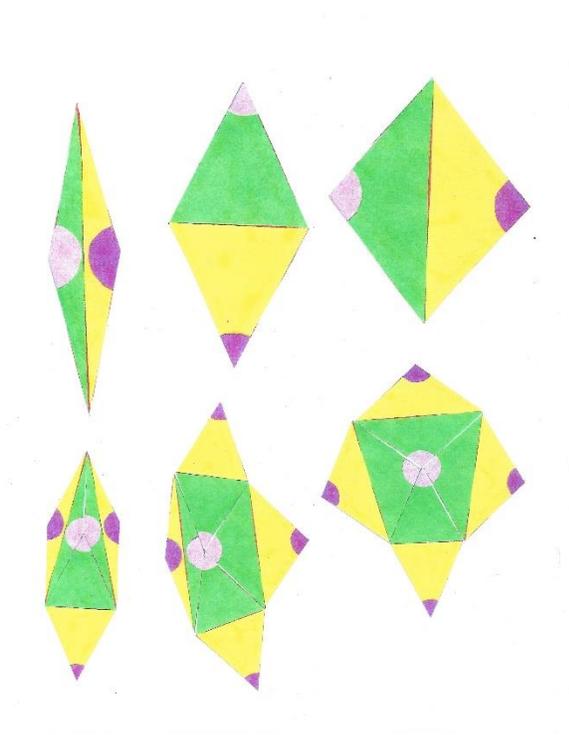

Fig. 14: A polka dot coloring of the minimal 7-fold tiling.

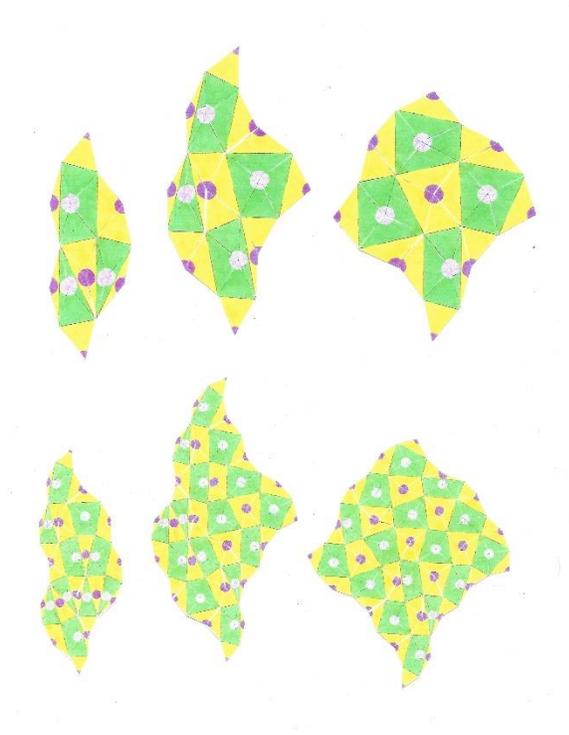

Fig. 15: The second and third generation of the polka dot variation.

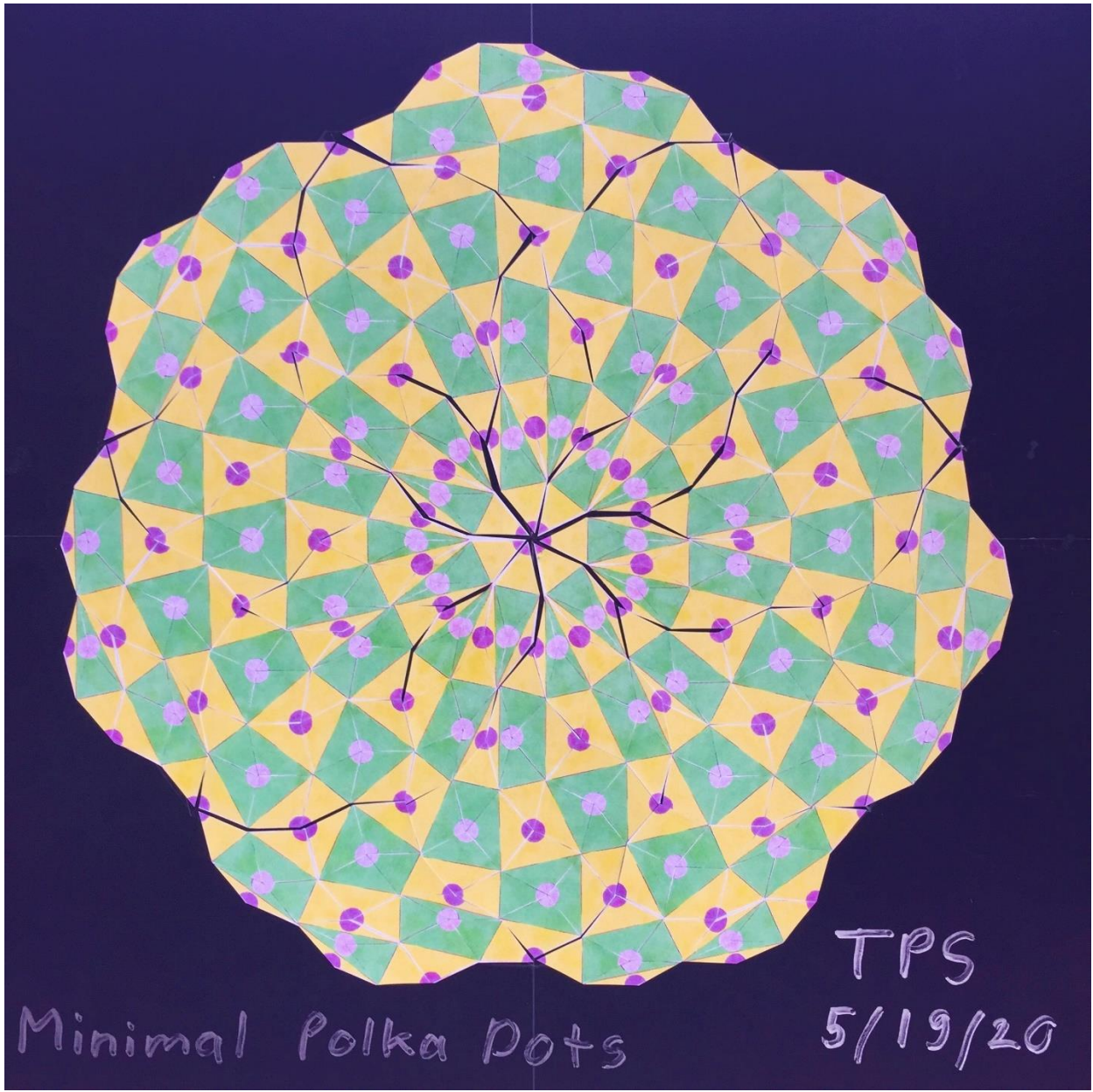

Fig 16: A minimal polka dot 7-fold mandala consisting of the three base rhombs.

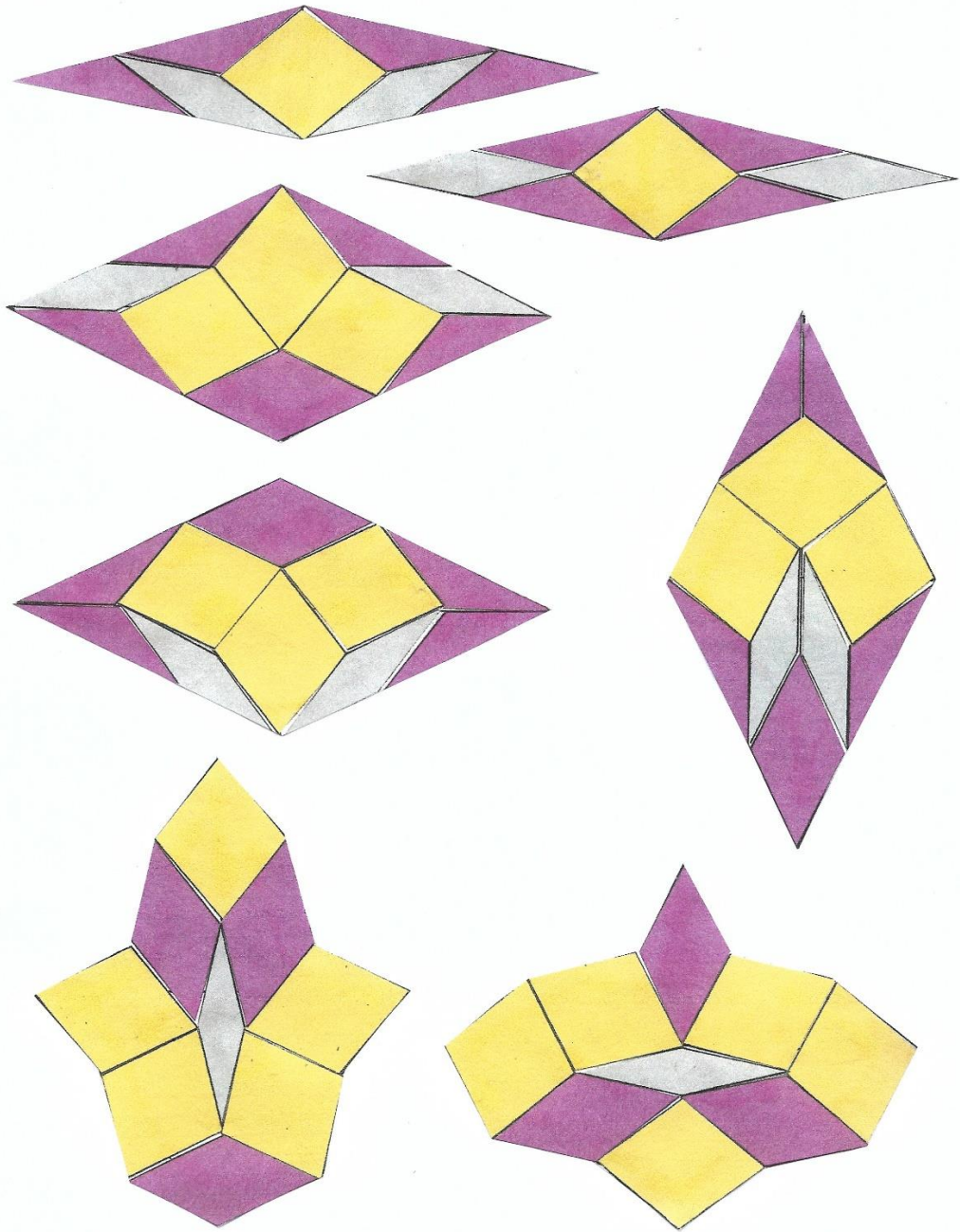

Fig. 17: Substitution rules for a scheme that uses 22 base tiles to reach the next stage. The new inflated edges have a length of δ = Φ + 1, where Φ is the long diagonal of the B tile. The inflation factor is δ = 2.802... However, there are no known matching rules that allow an inflation sequence that can be repeated.

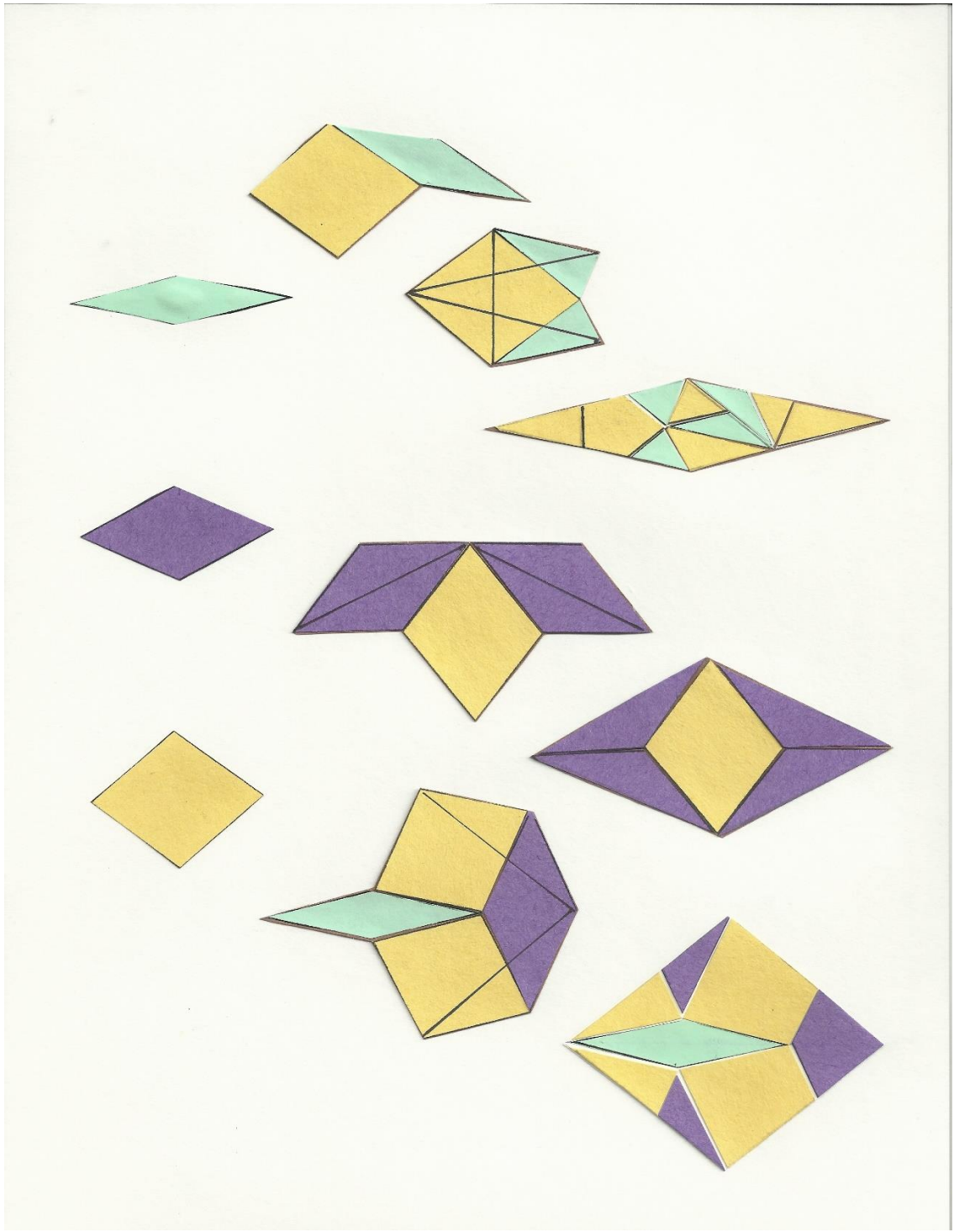

Fig. 18: The base tiles can fit into inflated tiles that are enlarged by the 7-fold magic number Φ with the help of some additional cuts. The inflation factor is δ = 1.802… This could be the smallest substitution set (9 tiles), but there are no known matching rules to repeat the inflation scheme.

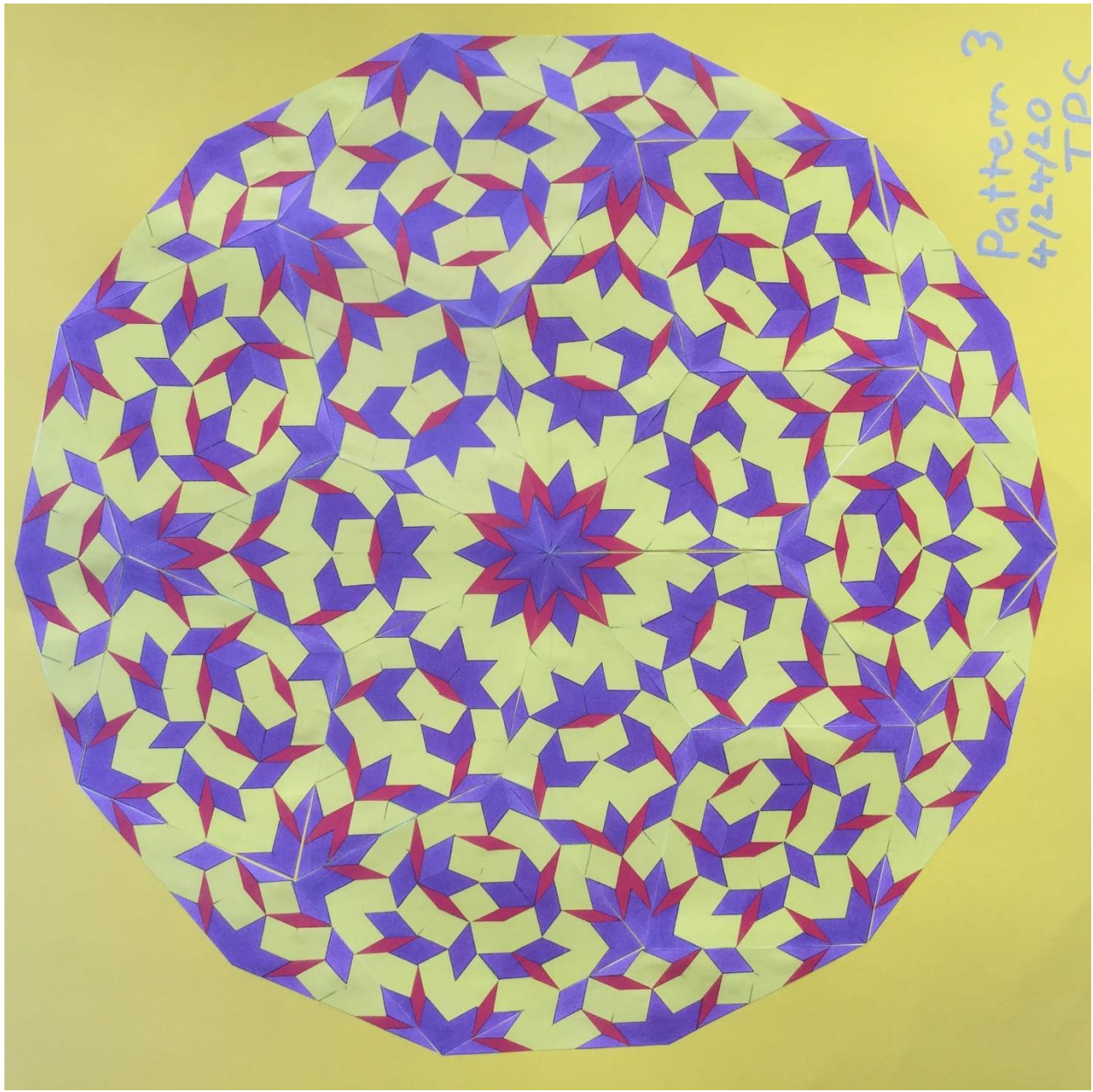

Fig. 19: A patch of an infinite 7-fold tiling nicknamed Sea Star and Crabs created during the coronavirus quarantine in a beach house on the Salish Sea.